\newif\ifcomments  
\commentsfalse
\documentclass{article}

\usepackage{amsmath,amssymb,pifont}
\usepackage{multicol}
\usepackage{amstext}
\usepackage{amsthm}
\usepackage{multirow}
\usepackage{booktabs}
\usepackage[skip=0pt]{subcaption}
\usepackage{times}
\usepackage{lipsum}
\usepackage[shortlabels]{enumitem}
\usepackage{cancel}
\usepackage{wrapfig}
\usepackage{array}
\usepackage{siunitx}
\usepackage{csvsimple}
\usepackage[multidot]{grffile}
\usepackage{bbm}
\usepackage{makecell}
\usepackage{bbm, dsfont}
\usepackage{mathtools}
\usepackage[dvipsnames]{xcolor}
\usepackage{comment}

\usepackage{geometry}
\usepackage{mathpazo}
\usepackage{enumitem}

\usepackage[multiple]{footmisc}
\usepackage{mathrsfs}
\usepackage{todonotes}
\usepackage{tikz}
\usepackage{ragged2e}   

\newcounter{myalgo}

\usepackage[ruled,vlined,linesnumbered]{algorithm2e}
\makeatletter
\let\c@algocf\c@myalgo
\makeatother
\usepackage{xfrac}

\usepackage{graphicx} %
\usepackage{color}

\usepackage{array}
\usepackage{amssymb}
\usepackage{amsmath}
\usepackage{xspace}
\usepackage{fancyhdr}
\usepackage{comment}
\usepackage{listings}
\usepackage[normalem]{ulem}
\usepackage{hyperref}
\usepackage{cleveref}

\usepackage[utf8]{inputenc}
\usepackage[table]{xcolor}
\usepackage{array}
\usepackage[most]{tcolorbox}
	
\hypersetup{final}

\newcommand{\calP}{\ensuremath{\mathcal{P}}}

\Crefname{lem}{Lemma}{Lemmas}
\Crefname{thm}{Theorem}{Theorems}
\Crefname{claim}{Claim}{Claims}
\Crefname{defn}{Definition}{Definitions}

\makeatletter
\newcommand{\vast}{\bBigg@{4}}
\newcommand{\Vast}{\bBigg@{5}}
\makeatother

\newcommand{\ex}[2]{{\ifx&#1& \mathbb{E} \else
\underset{#1}{\mathbb{E}} \fi \left[#2\right]}}
\newcommand{\pr}[2]{{\ifx&#1& \mathbb{P} \else
\underset{#1}{\mathbb{P}} \fi \left[#2\right]}}

 \newcommand{\ip}[2]{\left\langle #1, #2\right \rangle}

\usepackage{ulem}

\DeclarePairedDelimiterX{\infdivx}[2]{(}{)}{%
  #1\;\delimsize\|\;#2%
}

\newcommand{\mypar}[1]{\smallskip
	\noindent{\textbf{{#1}:}}}
	
\renewcommand{\epsilon}{\varepsilon}

\setlist{nolistsep}
\setlist[itemize]{noitemsep, topsep=0pt}

\setlist{nolistsep}
\setlist[itemize]{noitemsep, topsep=0pt}

\lstdefinestyle{mypython}{
  language=Python,
  backgroundcolor=\color{gray!10},
  basicstyle=\ttfamily\footnotesize,
  keywordstyle=\color{blue},
  commentstyle=\color{green!50!black},
  stringstyle=\color{red!60!black},
  showstringspaces=false,
  frame=single,
  breaklines=true
}

\usepackage{listings}
\usepackage{longtable}  

\newcommand{\ramsey}{\mathbf{R}}

\definecolor{light_blue}{RGB}{173, 216, 230} 
\definecolor{yellow_sota}{RGB}{255, 255, 0}   
\definecolor{light_green}{RGB}{144, 238, 144} 
\definecolor{light_red}{RGB}{240, 128, 128}   

\renewcommand{\themyalgo}{\thesection.\arabic{myalgo}}

\newtcolorbox[use counter=myalgo]{summaryboxgreen}[2][]{%
    colback=light_green!5,
    colframe=gray!60,
    coltitle=black,
    fonttitle=\bfseries,
    title={Algorithm~\themyalgo: #2}, 
    enhanced,
    attach boxed title to top left={yshift=-2mm, xshift=2mm},
    boxed title style={colback=light_green!20, colframe=gray!60, sharp corners},
    boxrule=0.5mm,
    arc=2mm,
    label={#1} 
}

\newtcolorbox[use counter=myalgo]{summaryboxblue}[2][]{%
    colback=light_blue!10,
    colframe=gray!60,
    coltitle=black,
    fonttitle=\bfseries,
    title={Algorithm~\themyalgo: #2}, 
    enhanced,
    attach boxed title to top left={yshift=-2mm, xshift=2mm},
    boxed title style={colback=light_blue!30, colframe=gray!60, sharp corners},
    boxrule=0.5mm,
    arc=2mm,
    label={#1} 
}

\newtcolorbox[use counter=myalgo]{summaryboxyellow}[2][]{%
    colback=yellow_sota!5,
    colframe=gray!60,
    coltitle=black,
    fonttitle=\bfseries,
    title={Algorithm~\themyalgo: #2}, 
    enhanced,
    attach boxed title to top left={yshift=-2mm, xshift=2mm},
    boxed title style={colback=yellow_sota!20, colframe=gray!60, sharp corners},
    boxrule=0.5mm,
    arc=2mm,
    label={#1} 
}

\newtcolorbox[use counter=myalgo]{summaryboxred}[2][]{%
    colback=light_red!5,
    colframe=gray!60,
    coltitle=black,
    fonttitle=\bfseries,
    title={Algorithm~\themyalgo: #2}, 
    enhanced,
    attach boxed title to top left={yshift=-2mm, xshift=2mm},
    boxed title style={colback=light_red!20, colframe=gray!60, sharp corners},
    boxrule=0.5mm,
    arc=2mm,
    label={#1} 
}

\crefname{myalgo}{Algorithm}{Algorithms}
\Crefname{myalgo}{Algorithm}{Algorithms}

\newcolumntype{L}[1]{>{\RaggedRight\arraybackslash\hspace{0pt}}p{#1}}
\definecolor{light_blue1}{RGB}{0, 100, 180}     
\definecolor{yellow_sota1}{RGB}{210, 180, 0}    
\definecolor{light_green1}{RGB}{0, 140, 0}      

\usepackage{fullpage}
\usepackage[numbers, sort, comma, square]{natbib}

\title{Reinforced Generation of Combinatorial Structures:\\ Ramsey Numbers}
\author{
Ansh Nagda\thanks{University of California, Berkeley, and Google DeepMind.}
\and
Prabhakar Raghavan\thanks{Google.}
\and
Abhradeep Thakurta\thanks{Google DeepMind.}
}
\date{}
 
\begin{document}

\maketitle
\thispagestyle{empty}
\begin{abstract}
We present improved lower bounds for nine classical Ramsey numbers: $\mathbf{R}(3, 13)$ is increased from $60$ to $61$, $\mathbf{R}(3, 18)$ from $99$ to $100$, $\mathbf{R}(4, 13)$ from $138$ to $139$, $\mathbf{R}(4, 14)$ from $147$ to $148$, $\mathbf{R}(4, 15)$ from $158$ to $159$, $\mathbf{R}(4, 16)$ from $170$ to $174$, $\mathbf{R}(4, 18)$ from $205$ to $209$, $\mathbf{R}(4, 19)$ from $213$ to $219$, and $\mathbf{R}(4, 20)$ from $234$ to $237$. These results were achieved using~\emph{AlphaEvolve}, an LLM-based code mutation agent. Beyond these new results, we successfully recovered lower bounds for all Ramsey numbers known to be exact, and matched the best known lower bounds across many other cases. These include bounds for which previous work does not detail the algorithms used.

Virtually all known Ramsey lower bounds for ``small'' numbers are derived computationally, with bespoke search algorithms each delivering a handful of results. AlphaEvolve is a single meta-algorithm yielding search algorithms for all of our results.
\end{abstract}
\clearpage
\newpage
\thispagestyle{empty}
\clearpage
\newpage
\pagenumbering{arabic}
\newpage


\section{Introduction}
\label{sec:intro}

Ramsey numbers have been extensively studied in the literature~\cite{radziszowski2024small}. In this work we focus on the classical graph formulation of Ramsey numbers, where $\ramsey(r,s)$ is the smallest number $n$ such that every undirected graph $G$ on $n$ vertices has either has a clique of size $r$, or an independent set of size $s$. While tight bounds on $\ramsey(r,s)$ are known for small values of $r$ and $s$ (namely, for $(3,s)$ with $3\leq s\leq 9$ and $(4,s)$ with $s\in\{4,5\}$~\cite{radziszowski2024small}), there is a large gap for many other values of $r$ and $s$. In this work, we focus on improving lower bounds on $\ramsey(r,s)$. If we can exhibit a graph with $n$ vertices that has no cliques of size $r$, nor any independent set of size $s$, we establish that $\ramsey(r,s)\geq n+1$.

Except for very small $r,s$ where analytical methods have worked, prior work has used a variety of search algorithms to computationally derive lower bounds on $\ramsey(r,s)$ (by generating graphs without the appropriate cliques and independent sets); some of these algorithms yielded the best known results for many pairs $r,s$~\cite{radziszowski2024small}.
We invoke AlphaEvolve~\cite{romera2024mathematical,novikov2025alphaevolve}, a code-mutation agent that uses LLMs (Large-Language-Models) to iteratively evolve code-snippets to generate better lower bounds. For each $r,s$, we use AlphaEvolve to generate search algorithms that produce larger and larger graphs without violating the clique and independence set constraints. (The exact experimental setup is a bit more elaborate.) Thus we use AlphaEvolve as a single meta-algorithm to generate search procedures for meeting or beating Ramsey lower bounds for many values of $r$ and $s$. Our work is similar in spirit to the discovery of extremal combinatorial objects in~\cite{nagda2025reinforced,georgiev2025mathematical}.


\subsection{Summary of Our Results}
\label{sec:results}
In this paper, we focus on~\emph{small} Ramsey numbers~\cite{radziszowski2024small} where $r,s\leq 22$. 
Our main results are that we improve $\ramsey(3, 13)$ from $60$~\cite{Kolodyazhny2015} to $61$, $\ramsey(3, 18)$ from $99$~\cite{Exoo2006} to $100$, $\ramsey(4,13)$ from $138$~\cite{exoo2015new} to $139$, $\ramsey(4, 14)$ from $147$~\cite{exoo2015new} to $148$,  $\ramsey(4, 15)$ from $158$~\cite{Tat2020}, $\ramsey(4, 16)$ from $170$~\cite{Tat2020} to $174$, $\ramsey(4, 18)$ from $205$~\cite{radziszowski2024small} to $209$, $\ramsey(4, 19)$ from $213$~\cite{radziszowski2024small} to $219$, and $\ramsey(4, 20)$ from $234$~\cite{radziszowski2024small} to $237$. Furthermore, we match the SoTA on many other cells in~\Cref{tbl:ramsey_lb}, including all cells for which the exact value of $\ramsey(r,s)$ is known. 
The search algorithms used to match or improve these 28 bounds can be categorized into four distinct families of initialization; these are summarized in Table~\ref{tab:ramsey_algos_grid} and detailed in the Appendix.
Note that for some of these cells, prior work only exhibits the graph implying the bound, without giving the search algorithm that discovered the graph~\cite{radziszowski2024small}.

\mypar{Note} Our constructions for the the cells $\ramsey(3, 13), \ramsey(3, 18), \ramsey(4,13), \ramsey(4,14)$, $\ramsey(4, 15)$, $\ramsey(4, 16)$, $\ramsey(4,18)$, $\ramsey(4,19)$, and $\ramsey(4,20)$ are made public on Github \href{https://github.com/google-research/google-research/tree/master/ramsey_number_bounds/improved_bounds}{\color{blue}\underline{at this link}}. Additionally, the search algorithms discovered by AI are also made public on Github \href{https://github.com/google-research/google-research/tree/master/ramsey_number_bounds/code}{\color{blue}\underline{at this link}}.

\begin{table}[h]
\centering
\renewcommand{\arraystretch}{1.5}
\setlength{\tabcolsep}{5pt}
\begin{tabular}{|c|c|c|c|c|c|c|c|c|c|c|c|c|c|c|c|c|c|c|c|c|}
\hline
\textbf{r \textbackslash s} & \textbf{3} & \textbf{4} & \textbf{5} & \textbf{6} & \textbf{7} & \textbf{8} & \textbf{9} & \textbf{10} & \textbf{11} & \textbf{12} & \textbf{13} & \textbf{14} & \textbf{15} & \textbf{16} & \textbf{17} & \textbf{18} & \textbf{19} & \textbf{20} & \textbf{21} & \textbf{22} \\ \hline
\textbf{3} & \cellcolor{light_blue}6 & \cellcolor{light_blue}9 & \cellcolor{light_blue}\hyperref[alg:r35]{14} & \cellcolor{light_blue}\hyperref[alg:r36]{18} & \cellcolor{light_blue}\hyperref[alg:r37]{23} & \cellcolor{light_blue}\hyperref[alg:r38]{28} & \cellcolor{light_blue}\hyperref[alg:r39]{36} & \cellcolor{yellow_sota}\hyperref[alg:r310]{40} & \cellcolor{yellow_sota}\hyperref[alg:r311]{47} &  
 & 
\cellcolor{light_green}\hyperref[alg:r313]{61} &  
 & 
& 
\cellcolor{yellow_sota}\hyperref[alg:r316]{82} & \cellcolor{yellow_sota}\hyperref[alg:r317]{92} & \cellcolor{light_green}\hyperref[alg:r318]{100} & \cellcolor{yellow_sota}\hyperref[alg:r319]{106} & \cellcolor{yellow_sota}\hyperref[alg:r320]{111} & \cellcolor{yellow_sota}\hyperref[alg:r321]{122} & \cellcolor{yellow_sota}\hyperref[alg:r322]{132} \\ \hline
\textbf{4} &  & \cellcolor{light_blue}\hyperref[alg:r44]{18} & \cellcolor{light_blue}\hyperref[alg:r45]{25} & \cellcolor{yellow_sota}\hyperref[alg:r46]{36} & \cellcolor{yellow_sota}\hyperref[alg:r47]{49} &  
 & 
 & \cellcolor{yellow_sota}\hyperref[alg:r410]{92} & & \cellcolor{yellow_sota}\hyperref[alg:r412]{128} & \cellcolor{light_green} \hyperref[alg:r413]{139} & \cellcolor{light_green}\hyperref[alg:r414]{148} & \cellcolor{light_green}\hyperref[alg:r415]{159} & \cellcolor{light_green}\hyperref[alg:r416]{174} & \cellcolor{yellow_sota}\hyperref[alg:r417]{200} & \cellcolor{light_green}\hyperref[alg:r418]{209} & \cellcolor{light_green}\hyperref[alg:r419]{219} & \cellcolor{light_green}\hyperref[alg:r420]{237} & &
 \\ \hline
\textbf{5} &  &  &  &  & \cellcolor{yellow_sota}\hyperref[alg:r57]{80} & \cellcolor{yellow_sota}\hyperref[alg:r58]{101} & \cellcolor{yellow_sota}\hyperref[alg:r59]{133} &  &  &  &  &  &  &  &  &  &  &  &  &  \\ \hline
\textbf{6} &  &  &  &  & \cellcolor{yellow_sota}\hyperref[alg:r67]{115} & &  &  &  &  &  &  &  &  &  &  &  &  &  &  \\ \hline
\end{tabular}
\caption{Ramsey Number $R(r, s)$ lower bounds achieved by AlphaEvolve. Highlights indicate performance status: \colorbox{light_blue}{\textbf{light blue}} for optimal lower bounds (matching upper bounds); \colorbox{yellow_sota}{\textbf{yellow}} for matching current best-known lower bounds (where the upper bound remains strictly higher); and \colorbox{light_green}{\textbf{green}} for new results that surpass the prior best lower bound. Algorithms are hyperlinked.}
\label{tbl:ramsey_lb}
\end{table}

\section{AlphaEvolve Meta-search, and Discovered Search Algorithms}
\label{sec:algo}
AlphaEvolve~\cite{novikov2025alphaevolve,romera2024mathematical} can be used to maintain a family of search algorithms (that search for graphs that establish Ramsey lower bounds), using an LLM to evolve them to find better search algorithms. Typically this process uses a synthetic objective function to guide AlphaEvolve, as a \textit{proxy} for the true objective (the size of the graph being discovered). This detail will become evident in our presentation of~\Cref{algo:ramseyAE} below. 

\mypar{AlphaEvolve experimental setup} \Cref{algo:ramseyAE} details the way this evolution develops search algorithms that try to improve the current SoTA for $\ramsey(r,s)$ (denoted as $n_{\texttt{SoTA}}$ in~\Cref{algo:ramseyAE}). Here are the core ideas:

\begin{itemize}
    \item \textbf{Evolution:} Maintain a population $\mathcal{P}$ of search algorithms. Use a performance-based $\texttt{Select}$ function to choose a parent algorithm from $\mathcal{P}$, then prompt an LLM to generate a mutated version $\texttt{p}_{\texttt{new}}$.
    \item \textbf{Execution:} Execute $\texttt{p}_{\texttt{new}}$ to get two graphs: a primary graph $G_1$ and a larger "prospect" graph $G_2$.
    \item \textbf{Scoring:} If $G_1$ is valid, it receives a score based on its size, heavily boosted if it is beyond $n_{\texttt{SoTA}}$. $G_2$ contributes an additive bonus inversely proportional to its violation count, guiding the search to push boundary limits. Formally, $\mathbb{E}_{viol}$ represents the expected total number of $r$-cliques and $s$-independent sets in a random graph of size $v_2$ (with edge probability 0.5). The score $S$ is updated with a term that increases linearly as the number of violations in $G_2$ decreases relative to the expected baseline $\mathbb{E}_{viol}$.
\end{itemize}
A central feature of~\Cref{algo:ramseyAE} is that it initializes with an empty baseline graph ($\texttt{p}_{\texttt{base}}$). AlphaEvolve then generates search programs $\texttt{p}_{\texttt{new}}$ aimed at discovering larger graphs. This approach encourages broad exploration, allowing the evolutionary process to learn how to modify and extend smaller, valid lower-bound constructions into larger ones. This is in contrast to an approach where we are specifically searching for graph constructions of a particular size~\cite{Exoo2006,exoo2015new,exoo2023lower}.

\begin{algorithm}[htb]
\caption{AlphaEvolve: Ramsey Program Search}
\KwIn{Ramsey parameters $(r, s, n_{\texttt{SoTA}})$, where $r$ represents the clique size, and $s$ represents the independent set size, and $n_{\texttt{SoTA}}$ is the current best graph size for $\ramsey(r,s)$ from~\cite{radziszowski2024small}}
\KwOut{Program $p^*$ with best graph}

\textbf{Init:} Population $\calP \gets \{\texttt{p}_{\texttt{base}}\}$ \tcp*[f]{$\texttt{p}_{\texttt{base}}$ returns empty graph}

\While{\text{compute budget remaining}}{
    $\texttt{p}_{\texttt{new}} \gets \texttt{LLM\_Mutation}(\texttt{Select}(\calP))$\\
    $(G_1, G_2)\gets \texttt{p}_{\texttt{new}}.\texttt{run}()$\\
    \uIf{$G_1$ has no $r$-cliques or $s$-independent sets}{
        $v_1, v_2 \gets |V(G_1)|, |V(G_2)|$\\
        $S \gets \begin{cases} 
            4\cdot v_1 & \text{if } v_1 > n_{\texttt{SoTA}} \\ 
            2\cdot v_1 & \text{if } v_1 = n_{\texttt{SoTA}} \\ 
            v_1 & \text{otherwise} 
        \end{cases}$  
        
        \BlankLine
        \tcp{Bonus: Reward low-violation "prospect" graphs}
        \If{$v_2 > v_1$}{
            $\mathbb{E}_{viol} \gets \frac{\binom{v_2}{r}}{2^{\binom{r}{2}}} + \frac{\binom{v_2}{s}}{2^{\binom{s}{2}}}$\\
            $S \gets S + \frac{1}{2}\max\left(0, 1 - \frac{\texttt{count\_viol}(G_2)}{\mathbb{E}_{viol}}\right)$\, where \texttt{count\_viol} returns the number of $r$-cliques and $s$-independent sets in the graph
        }
    }
    \Else{
        $S \gets -1$
    }
    $\texttt{score}(\texttt{p}_{\texttt{new}})\gets S$\\
    $\calP \gets \calP \cup \{(\texttt{p}_{\texttt{new}}, S)\}$
}
\KwRet{$\arg\max\limits_{\texttt{p} \in \calP} \texttt{score}(\texttt{p})$}
\label{algo:ramseyAE}
\end{algorithm}

Algorithms~\ref{alg:r313},~\ref{alg:r318},~\ref{alg:r413},~\ref{alg:r414},~\ref{alg:r415},~\ref{alg:r416},~\ref{alg:r418},~\ref{alg:r419}, and~\ref{alg:r420} detail the nine search algorithms that AlphaEvolve found for the cells $\ramsey(3,13)$, $\ramsey(3,18)$, $\ramsey(4,13)$, $\ramsey(4,14)$, and $\ramsey(4,15)$, $\ramsey(4,16)$, $\ramsey(4,18)$, $\ramsey(4,19)$ and $\ramsey(4,20)$ respectively. The detailed algorithms for all other colored cells in~\Cref{tbl:ramsey_lb} are presented in Algorithms~\ref{alg:r35}--\ref{alg:r67}. Note that each of the search algorithms seems specific to the corresponding cells; their success does not seem to transfer between cells, even when run for a reasonably long time. In~\Cref{tab:ramsey_algos_grid}, we summarize the initialization strategies used in the search algorithms.  We left out algorithms for $\ramsey(3,3)$ and $\ramsey(3,4)$ because AlphaEvolve just quoted the corresponding adjacency matrices from the LLM's memory. The search algorithms for each of the cells (while having some differences) are different variants of stochastic search.

If we take a closer look at the search algorithms discovered by AlphaEvolve for the colored cells in~\Cref{tbl:ramsey_lb}, a pattern starts emerging: a) For some of the cells,  a standard (stochastic) search from a generic starting point of Erd\H{o}s-R\'{e}nyi graphs suffices, b) For some, it was necessary to seed with specific algebraic graphs (e.g., Paley/cyclic), and c) For some the search had to be guided within some structural restrictions (e.g., staying within the family of cyclic graphs).
While clique and independent set counting is often a computational bottleneck in Ramsey constructions, the graph sizes for the instances in~\Cref{tbl:ramsey_lb} were relatively small. Thus, clique counting was not a limiting factor for these results.

\mypar{Comment on the use of AI summarization of the AlphaEvolve generated code} We utilized Gemini 3 to summarize the AlphaEvolve-generated code for all cells and to construct the classifications in~\Cref{tab:ramsey_algos_grid}. We scrutinized the underlying code to ensure that these AI-generated summaries reflect the search algorithms. In~\Cref{tab:ramsey_algos_grid} we see a bucketing of the various algorithms by Gemini, with the split being on the method by which each algorithm initialized its graph search.

\begin{table}[htb]
\centering
\small
\renewcommand{\arraystretch}{1.5}

\begin{tabular}{|>{\raggedright\arraybackslash}p{3.8cm}|>{\raggedright\arraybackslash}p{3.8cm}|>{\raggedright\arraybackslash}p{3.8cm}|>{\raggedright\arraybackslash}p{3.8cm}|}
\hline
\textbf{Random \& Stochastic Initialization} & \textbf{Paley \& Algebraic Seeding} & \textbf{Circulant \& Cyclic Bootstrap} & \textbf{Hybrid, Fractal \& Spectral Seeding} \\
\hline

\textbf{\textcolor{light_blue1}{R(3,5)}}, \textbf{\textcolor{light_blue1}{R(3,6)}}, 
\textbf{\textcolor{light_blue1}{R(3,7)}}, \textbf{\textcolor{light_blue1}{R(3,8)}},
\textbf{\textcolor{light_blue1}{R(4,5)}}
&
\textbf{\textcolor{light_blue1}{R(3,6)}}, \textbf{\textcolor{yellow_sota1}{R(3,10)}}, \textbf{\textcolor{yellow_sota1}{R(3,20)}}, \newline
\textbf{\textcolor{light_blue1}{R(4,4)}}, \textbf{\textcolor{yellow_sota1}{R(4,12)}}, 
\textbf{\textcolor{light_green1}{R(4,13)}}, \newline
\textbf{\textcolor{light_green1}{R(4,15)}},
\textbf{\textcolor{yellow_sota1}{R(6,7)}}
&
\textbf{\textcolor{light_blue1}{R(3,9)}}, \textbf{\textcolor{light_green1}{R(3,13)}}, 
\textbf{\textcolor{yellow_sota1}{R(3,16)}}, \newline
\textbf{\textcolor{yellow_sota1}{R(3,17)}}, 
\textbf{\textcolor{light_green1}{R(3,18)}}, \textbf{\textcolor{yellow_sota1}{R(3,19)}}, \newline
\textbf{\textcolor{yellow_sota1}{R(3,21)}}, \textbf{\textcolor{yellow_sota1}{R(3,22)}}, 
\textbf{\textcolor{yellow_sota1}{R(4,6)}}, \newline\textbf{\textcolor{yellow_sota1}{R(4,7)}}, 
\textbf{\textcolor{yellow_sota1}{R(4,10)}}, \textbf{\textcolor{light_green1}{R(4,14)}},
\newline
\textbf{\textcolor{light_green1}{R(4,16)}}, \textbf{\textcolor{yellow_sota1}{R(4,17)}}, \textbf{\textcolor{light_green1}{R(4,18)}},\newline \textbf{\textcolor{light_green1}{R(4,19)}}, \textbf{\textcolor{light_green1}{R(4,20)}}, \textbf{\textcolor{yellow_sota1}{R(5,7)}},\newline
\textbf{\textcolor{yellow_sota1}{R(5,8)}}, \textbf{\textcolor{yellow_sota1}{R(5,9)}}
&
\textbf{\textcolor{light_blue1}{R(3,6)}}, \textbf{\textcolor{yellow_sota1}{R(3,11)}}
\\ \hline\hline

\begin{itemize}[leftmargin=*, nosep, label=\tiny\textbullet]
    \item \textbf{Initialization:} Random graphs ($G(n, p)$) or empty/greedy baseline.
\end{itemize}
&
\begin{itemize}[leftmargin=*, nosep, label=\tiny\textbullet]
    \item \textbf{Initialization:} Explicit seeding with Paley graphs, cubic, and quadratic residue graphs.
\end{itemize}
&
\begin{itemize}[leftmargin=*, nosep, label=\tiny\textbullet]
    \item \textbf{Initialization:} Bootstrapped from circulant graphs, sum-free sets, or cyclic constructions.
\end{itemize}
&
\begin{itemize}[leftmargin=*, nosep, label=\tiny\textbullet]
    \item \textbf{Initialization:} Complex mixtures: Fractal/self-similar construction, spectral properties, or hybrid pools.
\end{itemize}
\\ \hline

\end{tabular}
\caption{Taxonomy of Ramsey Search Algorithms by Initialization Method}
\label{tab:ramsey_algos_grid}
\end{table}



\section{Comparison to Prior Work}
\label{sec:comp}

We compare our results against prior work that established the previous state-of-the-art (SoTA) for the entries in~\Cref{tbl:ramsey_lb}. Our primary reference for these baselines is~\cite{radziszowski2024small}. Since many historical results rely on personal communications without published algorithmic details, we restrict our comparison to cases where: a) a public reference is available, and b) a computational approach was employed to achieve the lower bound.

\begin{itemize}
    \item $\ramsey(4,13), \ramsey(4,14),\ramsey(4,15)$,$\ramsey(4,15)$, $\ramsey(4,16)$, $\ramsey(4,18)$ and $\ramsey(6,7)$: Except for $\ramsey(4,s), s\in\{15,16,18\}$, the previous SoTA for these cells were established in~\cite{exoo2015new}. For $\ramsey(4,15)$ and $\ramsey(4,16)$ the prior SoTA was established in~\cite{Tat2020}.  We improved the bound for $\ramsey(4,13)$ (and $\ramsey(4,14)$), while matching the bound for $\ramsey(6,7)$. Notably, the search algorithms discovered by AlphaEvolve for $\ramsey(4,13)$ and $\ramsey(6,7)$ (Algorithms~\ref{alg:r413} and~\ref{alg:r67}) initialize with the same algebraic structures—cubic residue and Paley graphs, respectively—as those used in~\cite{exoo2015new}. For $\ramsey(4,14)$, Exoo~\cite{exoo2015new} initialized the search with a cubic graph, whereas AlphaEvolve (Algorithm~\ref{alg:r414}) used a class of circulant graphs with varying periodicities. 
    The subsequent optimization strategies for all the cells differ significantly: AlphaEvolve utilizes evolved heuristics rather than the standard simulated annealing employed in~\cite{exoo2015new}. For $\ramsey(4,18)$, the previous SoTA~\cite{radziszowski2024small} was obtained via differencing between other Ramsey cells. However, AlphaEvolve developed a new ground-up heuristic for this cell.
    
    Our approach for $\ramsey(4,15)$ (\Cref{alg:r415}) diverges from the Cayley graph initialization of~\cite{exoo2015new} (which is the SoTA prior to~\cite{Tat2020}), starting instead from a Paley graph. We unfortunately could not compare to~\cite{Tat2020} as it is referred to as personal communication in~\cite{radziszowski2024small}.

    \item $\ramsey(4,6)$: We match the SoTA for this cell, which was achieved by~\cite{Exoo2012}. While~\cite{Exoo2012} initializes the search with a random graph, AlphaEvolve opts for a cyclic graph initialization. Despite this difference, the subsequent search procedure discovered by AlphaEvolve (\Cref{alg:r46}) bears structural similarities to the heuristic search described in~\cite{Exoo2012}.

    \item $\ramsey(4,10)$: The current SoTA for this cell is attributed to~\cite{Harborth2003}, who employed a branch-and-bound search on circulant colorings. In contrast, the algorithm discovered by AlphaEvolve (\Cref{alg:r410}) initiates by sampling from a distribution over circulant graphs, followed by a custom local search that integrates sophisticated tabu search mechanisms with sequential growth.
\end{itemize}

This comparison demonstrates that for many cells --- even those where AlphaEvolve only matched the SoTA --- AlphaEvolve's search algorithms differ materially from prior methods, despite often sharing initialization families (e.g., algebraic or cyclic graphs). At a meta-level, the algorithms generated by AlphaEvolve exhibit three salient properties: a) they exploit known algebraic structures relevant to the specific cell; b) unlike many human-designed algorithms, they typically chain multiple heuristics or execute them in parallel; and c) they incorporate custom heuristics for approximate clique and independent set counting to accelerate candidate evaluation. Furthermore, for cases such as $\ramsey(4,10)$ (see~\Cref{alg:r410}), AlphaEvolve's search strategies appear novel and, to the best of our knowledge, are absent from existing literature.

\mypar{Prior work on AI and extremal combinatorics} Recent work has studied the use of AI in extremal combinatorics: {finding counterexamples to graph theory conjectures}~\cite{wagner2021constructions}, {the capset problem}~\cite{romera2024mathematical}, {constructing maximum grid subsets avoiding isosceles triangles}~\cite{charton2024patternboost}, {discovering extremal Ramanujan graphs}~\cite{nagda2025reinforced}, and {improving bounds for the finite field Kakeya problem}~\cite{georgiev2025mathematical}. Apart from~\cite{wagner2021constructions}, all the the other results (including ours) have used some variant of AlphaEvolve to discover algorithms that eventually search for extremal combinatorial structures. This contrasts with some other advancements in the space of AI and math/theoretical computer science~\cite{Bubeck2025_GPT5_Proof,woodruff2026accelerating} where new theorems were discovered by directly prompting LLMs (large-language-models). It is unclear at this point whether directly prompting an LLM will result in discovering extremal combinatorial structures.
\section{Conclusion}
\label{sec:conclusion}

LLM-based systems like AlphaEvolve have been applied to a variety of problems in mathematics and computer science (see~\cite{nagda2025reinforced} for a detailed survey). The area of using AI for Math/CS research is a fast evolving landscape. We approach it from the point of view of making progress on classic and well-studied problems, to derive results that stand the test of time. 

Our work lies within the domain of \emph{extremal combinatorics} and is most closely related to~\cite{nagda2025reinforced}, which generated gadgets and extremal graphs to establish new hardness of approximation results. However, while~\cite{nagda2025reinforced} required significant human effort to design proof structures amenable to AlphaEvolve, here the human effort is in crafting the overarching \Cref{algo:ramseyAE}. Unlike prior computational attempts at improving Ramsey lower bounds~\cite{radziszowski2024small} --- which typically rely on executing bespoke, human-designed heuristics --- AlphaEvolve automates the heuristic discovery process: it searches for \emph{novel search algorithms}.

In this work, we focus on Ramsey lower bounds, where the objective is to generate a witness graph of a target size that avoids specific cliques and independent sets. The complementary challenge of improving Ramsey upper bounds requires a fundamentally different approach: demonstrating the \emph{non-existence} of any such graph larger than a specific threshold. Unlike lower bound searches, this cannot be solved by constructing a single example. Recently, formal methods have been successfully employed to prove that $\ramsey(4,5)=25$~\cite{gauthier2024formal}.



\bibliographystyle{alpha}
\bibliography{reference}
\appendix
\begin{summaryboxgreen}[alg:r313]{Search algorithm for {$\ramsey(3,13)$} lower bound}

\begin{itemize}
\item \textbf{Phase 1: Cyclic Bootstrap:}The search begins by rapidly identifying a large base graph $G_{\texttt{start}}$ using a randomized cyclic construction. Edges are defined by a difference set $S$, where $(u, v) \in E \iff |u-v| \pmod n \in S$. The algorithm maximizes $n$ such that the graph remains $(3, 13)$-free.

\item \textbf{Phase 2: Iterative Expansion:} 
Starting from $G_{\texttt{curr}} \leftarrow G_{\texttt{start}}$, the algorithm incrementally increases the graph size to $n+1$. A new vertex $v_{\texttt{new}}$ is added, and its connections are optimized to maintain the $(3, 13)$-free property using advanced heuristics.

\item \textbf{Initialization (Hitting Sets):} 
The neighborhood of $v_{\texttt{new}}$ is explicitly constructed to "hit" (intersect) existing independent sets of size $s-1$ (12-anticliques). By ensuring $v_{\texttt{new}}$ connects to at least one vertex in every such set, the algorithm prevents them from growing to size $s=13$:
$$N(v_{\texttt{new}}) \leftarrow \{ u \in V \mid u \text{ intersects } C \in \mathcal{C}_{12}(G_{\texttt{curr}}) \}$$

\item \textbf{Candidate Selection:} 
Before full optimization, multiple initialization strategies (Hitting Sets, Max Independent Set, Random) generate candidate neighborhoods. Each candidate undergoes a short ``Mini-SA'' burst—a truncated annealing run with an accelerated cooling schedule—to empirically evaluate its quality and discard unpromising seeds before the expensive main optimization begins. The candidate with the lowest energy is selected for the main refinement phase.

\item \textbf{Adaptive Simulated Annealing (SA):} 
The selected graph is refined to minimize a weighted energy function:
$$E(G) = \lambda_{K_3} \cdot \text{count}(K_3) + \lambda_{I_{13}} \cdot \text{count}(I_{13})$$
The weights $\lambda$ are \textbf{adaptive}: if $\text{count}(K_3)$ dominates, $\lambda_{K_3}$ increases dynamically to force the solver to prioritize destroying triangles.

\end{itemize}
\end{summaryboxgreen}

\begin{summaryboxgreen}[alg:r318]{Search algorithm for {$\ramsey(3,18)$} lower bound}
\begin{itemize}
    \item \textbf{Main Loop (Iterative Deepening):} Initialize target size $n \gets \max(60, |V(G_{\texttt{best}})| + 1)$. The algorithm iterates through increasing $n$, attempting to construct a valid $(3, 18)$-free graph using three strategies.

    \item \textbf{Strategy 1: Cyclic Construction:} 
    Generates a maximal sum-free set $S \subset \{1, \dots, \lfloor n/2 \rfloor\}$. The graph is constructed such that $(u, v) \in E \iff |u-v| \pmod n \in S$. This guarantees $K_3$-freeness ($r$-free).

    \item \textbf{Strategy 2: Block Construction:} 
    Decomposes $n$ into $n_b$ blocks of size $k$. Uses two sum-free sets: $S_{\texttt{intra}}$ for edges within blocks and $S_{\texttt{inter}}$ for edges between blocks $B_i$ and $B_j$.
    $$ E_{\texttt{block}} = E(S_{\texttt{intra}}) \cup E(S_{\texttt{inter}}) $$

    \item \textbf{Strategy 3: Adaptive Growth:} 
    Initializes with the best valid graph found so far. Adds a new vertex $v_{\texttt{new}}$ and connects it to existing vertices $u \in V$ greedily, provided no triangle is formed:
    $$ N(v_{\texttt{new}}) \leftarrow \{ u \mid \nexists w \in N(v_{\texttt{new}}) \text{ s.t. } (u, w) \in E \} $$
    
    \item \textbf{Filtration \& Validation:}
    Candidates are first filtered by a heuristic $\texttt{estimate}\  \alpha(G)$. If the heuristic fails to find an independent set of size $s$, an exact solver computes $\alpha(G)$. If $\alpha(G) < s$, update $G_{\texttt{best}} \leftarrow G$ and increment $n$.

\end{itemize}
\end{summaryboxgreen}

\begin{summaryboxgreen}[alg:r413]{Search algorithm for {$R(4,13)$}}
    \begin{itemize}
        \item \textbf{Algebraic Bootstrap ($C_{127}^3$):}
        The search is initialized with a highly specific algebraic construction: the Cayley graph formed by cubic residues in the finite field $\mathbb{F}_{127}$.
        $$ V = \mathbb{Z}_{127}, \quad (u,v) \in E \iff (u - v) \pmod{127} \in \{x^3 \mid x \in \mathbb{Z}_{127}^*\} $$
        This provides a starting graph of size 127 that is likely very close to valid.

        \item \textbf{Incremental Violation Tracking:}
        Unlike previous approaches that only count defects, this algorithm maintains \textbf{explicit sets} of all current 4-cliques and 13-anticliques (e.g., `Set[FrozenSet[int]]`). This allows for precise, incremental updates of the defect lists upon every edge flip without full graph re-scanning.

        \item \textbf{Two-Stage Move Evaluation:}
        To handle the large neighborhood of potential moves, the algorithm uses a filter:
        \begin{enumerate}
            \item \textbf{Heuristic Filter:} Scores candidates using a fast proxy (counting common neighbors for $\Delta K_4$ and common non-neighbors for $\Delta I_{13}$ risk).
            \item \textbf{Exact Delta:} Performs the expensive subgraph isomorphism check only on the top survivors of the heuristic filter.
        \end{enumerate}

        \item \textbf{Perturbation-Based Extension:}
        When extending the graph ($n \to n+1$), the new vertex vector is generated by \textbf{cloning and perturbing} an existing row (flipping $\approx 5\text{-}15\%$ of edges). This assumes the optimal connection pattern for a new node is likely similar to valid patterns already present in the graph. The score for a candidate extended graph is computed quickly by making use of prior knowledge of the violations in the base graph.

        \item \textbf{Strategic Kick (Escape):}
        If the Tabu search stagnates (no improvement for 500 iterations), a "Strategic Kick" is applied: the algorithm forces a short sequence of random edge flips (ignoring objective value) to jump out of the local basin of attraction.
    \end{itemize}
\end{summaryboxgreen}

\begin{summaryboxgreen}[alg:r414]{Search algorithm for {$\ramsey(4,14)$} lower bound}\begin{itemize}\item \textbf{Orbit-Based Circulant Construction:} The search is restricted to highly symmetric, vertex-transitive circulant graphs. Instead of searching the full graph space, the algorithm generates a difference set $S \subset \mathbb{Z}_n$ that fully defines the graph: $(u, v) \in E \iff (u-v) \pmod n \in S$. The set $S$ is constructed as a union of orbits under modular multiplication by a generator $g$, including additive inverses to ensure the graph is undirected.\item \textbf{Vertex-Transitive Reduction:} Because the graph is vertex-transitive, global structural constraints are reduced to localized checks around a single vertex (vertex 0). Here, $S$ acts exactly as the neighborhood of vertex 0:
\begin{align*}K_4 \text{ freeness} &\iff \text{No } K_3 \text{ in the subgraph induced by the neighborhood } S \\
I_{14} \text{ freeness} &\iff \text{No } I_{13} \text{ in the subgraph induced by the non-neighborhood } \mathbb{Z}_n \setminus (S \cup {0})\end{align*}
\item \textbf{Stochastic Orbit Sampling:} When the number of generated orbits $m$ is small ($m \le 14$), the algorithm exhaustively evaluates all $2^m$ combinations to form the connection set $S$. For larger orbit spaces, it stochastically samples unions of orbits that yield a target neighborhood degree (restricting the size of $S$ to be between 24 and 48), rapidly skipping configurations outside this viable density window.
\item \textbf{High-Performance Bitset Filtration:} The localized $K_3$ and $I_{13}$ subgraph isomorphism checks are aggressively accelerated using custom bitwise operations. Adjacency lists for the subgraphs are compressed into bitmasks, allowing the recursive clique searches to use fast bit-shifts and masking to evaluate thousands of orbit combinations per second.\end{itemize}\end{summaryboxgreen}

\begin{summaryboxgreen}[alg:r415]{Search algorithm for {$\ramsey(4,15)$} lower bound}
    \begin{itemize}
        \item \textbf{Harmonic Genetic Memory:}
        The algorithm maintains a global ``gene pool'' across generations that records the frequency of successful edges and ``orbits'' (cyclic distances).
        $$ \mathcal{M}_{\texttt{edge}}[u,v] \gets \mathcal{M}_{\texttt{edge}}[u,v] + 1 \quad \text{if } (u,v) \in E(G_{\texttt{valid}}) $$
        $$ \mathcal{M}_{\texttt{orbit}}[d] \gets \mathcal{M}_{\texttt{orbit}}[d] + 1 \quad \text{if } \exists (u,v) \in E(G_{\texttt{valid}}) \text{ s.t. } |u-v| = d $$
        
        \item \textbf{Spectral Initialization:}
        New candidates for size $n$ are initialized either via algebraic constructions (Generalized Paley/Quadratic Residue graphs on primes $p \approx n$) or by extending $G_{\texttt{best}}$ using probabilities derived from the harmonic memory $\mathcal{M}$.
        
        \item \textbf{Tabu-Enhanced Simulated Annealing:}
        The graph undergoes local search to minimize strictly the clique count $K_4$. Moves (edge flips) are managed by a Tabu list to prevent cycling, with an adaptive cooling schedule that slows down when $K_4$ counts approach zero.
        
        \item \textbf{Look-Ahead \& Harmonic Scoring:}
        When removing an edge to break a $K_4$, candidates are scored not just by $\Delta K_4$, but by a ``Look-Ahead'' penalty (risk of creating large independent sets) and a ``Harmony'' bonus (preserving edges/orbits with high history in $\mathcal{M}$):
        $$ \text{Score}(e) = \Delta K_4 - \lambda_1 \cdot \text{Risk}(I_{\texttt{neighbor}}) - \lambda_2 \cdot \mathcal{M}(e) $$
        
        \item \textbf{Harmonic Tunneling:}
        If the search stagnates with $K_4 > 0$, the algorithm identifies a ``toxic orbit'' $d^*$—the cyclic distance most frequent within the remaining cliques. It then performs a massive mutation by flipping \textit{all} edges $(u,v)$ where $|u-v| = d^*$, tunneling out of the local minimum.
        
        \item \textbf{Constraint Toggle Strategy:}
        The search prioritizes $K_4=0$. Once satisfied, it switches to breaking Independent Sets ($I_{15}$). It identifies a violating set $S$ via heuristic sampling and adds an edge $(u,v)$ within $S$ that minimizes the creation of new $K_4$s.
    \end{itemize}
\end{summaryboxgreen}

\begin{summaryboxgreen}[alg:r416]{Search algorithm for {$\ramsey(4,16)$} lower bound}
\begin{itemize}
\item \textbf{Circulant Graph Representation:} The search is restricted to highly symmetric, vertex-transitive circulant graphs. A graph of size $n$ is fully defined by a difference set $S \subset \{1, \dots, \lfloor n/2 \rfloor\}$, where edges are formed if the modular distance between vertices belongs to $S$. The algorithm incrementally tests graph sizes $n \in [155, 180]$, carrying over successful difference sets to seed the next iteration.
\item \textbf{Vertex-Transitive Reduction:} Due to the graph's symmetry, global structural constraints are simplified to localized checks strictly around vertex 0. The presence of forbidden subgraphs is mapped directly to the neighborhood and non-neighborhood of a single vertex:
\begin{align*}
K_4 \text{ freeness} &\iff \text{No } K_3 \text{ in the subgraph induced by the neighborhood } S \\
I_{16} \text{ freeness} &\iff \text{No } I_{15} \text{ in the subgraph induced by the non-neighborhood } \mathbb{Z}_n \setminus (S \cup \{0\})
\end{align*}
\item \textbf{Simulated Annealing with Guided Mutations:} The connection set $S$ is optimized using simulated annealing, minimizing an energy function based on the count of $K_4$ subgraphs or the maximum independent set size. To accelerate convergence, mutations are aggressively guided (with 75\% probability) to directly target and disrupt specific violating structures, such as modifying the distance between vertices identified in a $K_4$ generator or an $I_{16}$ witness.
\item \textbf{Automorphism Teleportation and Bitset Filtration:} The localized subgraph checks are heavily optimized using custom bitwise operations for rapid clique and anticlique backtracking. To escape deep local minima after prolonged stagnation, the algorithm applies "symmetry-invariant automorphism teleportation," mapping the difference set $S \to (S \cdot k) \pmod n$ (where $\gcd(k, n) = 1$) to jump to an isomorphic graph state and resume searching.
\end{itemize}
\end{summaryboxgreen}

\begin{summaryboxgreen}[alg:r418]{Search algorithm for {$\ramsey(4,18)$} lower bound}
\begin{itemize}
\item \textbf{Circulant Core Generation via Orbits:} The initial search explores highly symmetric, vertex-transitive circulant graphs defined by a difference set $S$. To drastically compress the search space, distances are grouped into orbits under multiplicative subgroup actions, allowing the algorithm to construct graphs by toggling entire symmetry orbits rather than individual edges.
\item \textbf{Multi-Phase Simulated Annealing:} The circulant optimization is split into a "Rapid Scan" across graph sizes $N \in [195, 225]$ and a subsequent "Deep Dive" focusing on the most promising candidates. The annealing engine uses specialized moves like "Harmonic Orbit Translocation" to swap complementary orbit states, evaluated via high-speed bitwise triangle and clique counters.
\item \textbf{Non-Circulant Metamorphosis (Ghost Vertex Injection):} After isolating the maximal valid circulant cores, the algorithm intentionally breaks global symmetry to extend the graph further. It employs an SMT-inspired local search to iteratively attach "ghost vertices," dynamically flipping the new vertex's adjacencies to minimize induced $K_4$ and $I_{18}$ structural violations.
\item \textbf{Beam Search and Subgraph Caching:} The asymmetric extension phase is orchestrated using a beam search, maintaining a pool of the deepest valid parent graphs. The edge-assignment solver is heavily optimized by caching all existing $K_3$ and $I_{17}$ subgraphs as bitmasks, allowing it to evaluate and resolve constraint violations at a massive scale during the metamorphosis step.
\end{itemize}
\end{summaryboxgreen}

\begin{summaryboxgreen}[alg:r419]{Search algorithm for {$\ramsey(4,19)$} lower bound}
\begin{itemize}
\item \textbf{Dynamic Graph Sizing:} The search explores symmetric circulant graphs defined by a difference set $S$. It probabilistically jumps between graph sizes around $n \approx 216$, caching promising states and proportionally scaling successful difference sets ($S_{new} \approx S_{best} \cdot n_{new} / n_{best}$) to seed new sizes.
\item \textbf{Two-Tiered Evaluation:} Global constraints are reduced to fast bitwise neighborhood checks. The energy function heavily penalizes $K_4$ subgraphs. Only strictly $K_4$-free graphs trigger a rigorous, dynamically scaled backtracking search on the complement graph to verify $I_{19}$ freeness.
\item \textbf{Simulated Annealing \& Teleportation:} Optimization uses simulated annealing with an adaptive cooling rate that decelerates during plateaus. To escape deep local minima, a "Quantum Teleportation" mutation periodically replaces two distances $d_1, d_2 \in S$ with their "entangled" sum and difference.
\item \textbf{Guided Disruption:} Mutations explicitly target structural violations by dropping edges from $K_4$ generators or injecting edges to break $I_{19}$ cliques (verified via look-ahead checks). Periodic automorphism jumps $S \to (S \cdot k) \pmod n$ shift the search to isomorphic states to prevent stagnation.
\end{itemize}
\end{summaryboxgreen}

\begin{summaryboxgreen}[alg:r420]{Search algorithm for {$\ramsey(4,20)$} lower bound}
\begin{itemize}
\item \textbf{Prioritized Search \& Warm-Starting:} The algorithm targets a specific priority list of graph sizes ($n \in [180, 240]$). It accelerates convergence using "warm starts," seeding new searches by proportionally scaling the best known difference set from a previously explored size.
\item \textbf{Hierarchical Energy Evaluation:} The objective function $E = 1000 \cdot |K_4| + |I_{max}|$ strictly prioritizes $K_4$ freeness. Fast bitwise intersections count $K_4$ subgraphs; if manageable, an exact Bron-Kerbosch backtracking solver evaluates the maximum independent set $I_{max}$ on the complement graph.
\item \textbf{Advanced Mathematical Mutations:} Alongside random tweaks, the engine uses two advanced operators:
\begin{itemize}
\item \textit{Spectral Mutation:} Evaluates the graph's eigenvalues via a cosine transform, probabilistically adding or removing edges based on their gradient impact to break cliques or independent sets.
\item \textit{Harmonic Mutation:} Toggles distances along specific arithmetic progressions to disrupt periodic geometric structures.
\end{itemize}
\item \textbf{Targeted Disruption \& Extension Phase:} Mutations force localized edge flips directly on vertices involved in $K_4$ or $I_{20}$ violations. Once a valid "G1" core graph is discovered, a "G2" phase immediately attempts to extend the lower bound to sizes $n+1$ through $n+4$ using the optimized core as a starting seed.
\end{itemize}
\end{summaryboxgreen}

\begin{summaryboxblue}[alg:r35]{Search algorithm for {$R(3,5)$} lower bound}
    \begin{itemize}
        \item \textbf{Sequential Stochastic Search:}
        The algorithm iterates through increasing graph sizes $n$, starting from $n=10$. For each size, it initializes a random graph $G \sim G(n, 0.5)$ and attempts to refine it into a valid $(3, 5)$-free graph.

        \item \textbf{Objective Function:}
        The minimization target is the unweighted sum of forbidden subgraphs. The energy of a state is calculated by a full traversal of the graph:
        $$ E(G) = \text{count}(K_3) + \text{count}(I_5) $$

        \item \textbf{Simulated Annealing:}
        The graph evolves via a standard annealing process. A candidate move consists of flipping a single random edge $(u, v)$. Moves are accepted based on the Metropolis criterion with geometric cooling ($T_{k+1} = \alpha \cdot T_k$).

        \item \textbf{Termination \& Progression:}
        If the energy reaches $E(G)=0$, the graph is saved as the new $G_{\texttt{best}}$, and the target size is incremented ($n \gets n+1$). If the cooling schedule completes with $E(G) > 0$, the sequential search terminates.
    \end{itemize}
\end{summaryboxblue}

\begin{summaryboxblue}[alg:r36]{Search algorithm for {$R(3,6)$}}
    \begin{itemize}
        \item \textbf{Hybrid Initialization Strategy:}
        For each target size $n$, the algorithm selects from three initialization methods:
        \begin{enumerate}
            \item \textit{Constructive:} Extends the previous best graph $G_{n-1}$ by adding a vertex with sparse random connections ($p \approx 0.45$).
            \item \textit{Algebraic:} Constructs a Paley graph if $n$ is a prime power and $n \equiv 1 \pmod 4$. 
            \item \textit{Stochastic:} Generates random graphs $G(n, p)$ with varying densities ($p \in \{0.4, 0.5, 0.6\}$).
        \end{enumerate}

        \item \textbf{Asymmetric Cost Function:}
        The Simulated Annealing process minimizes a weighted energy function. The penalty for 3-cliques is set $100\times$ higher than for independent sets, forcing the solver to prioritize the strict $K_3$-free constraint:
        $$ E(G) = 100 \cdot \text{count}(K_3) + \text{count}(I_6) $$

        \item \textbf{Sequential Growth:}
        The search operates sequentially from $n=12$ to $n=18$. The algorithm iterates to size $n+1$ only after successfully finding a valid zero-defect graph for size $n$.

        \item \textbf{Simulated Annealing Refinement:}
        Refinement is performed via standard SA with geometric cooling ($T_{\texttt{new}} = 0.9999 \cdot T_{\texttt{curr}}$). The high triangle penalty effectively restricts the search space to ``nearly'' triangle-free graphs.
    \end{itemize}
\end{summaryboxblue}

\begin{summaryboxblue}[alg:r37]{Search algorithm for {$R(3,7)$}}
    \begin{itemize}
        \item \textbf{Iterative Deepening:}
        The search begins at $n=22$ (a likely lower bound for efficient search) and increments $n$ only after a valid $(3, 7)$-free graph is found.

        \item \textbf{Incremental Initialization:}
        For the initial size ($n=22$), the graph is generated randomly ($p=0.5$). For subsequent sizes ($n > 22$), the algorithm initializes by extending the previous best graph $G_{n-1}$, adding a new vertex with random connections ($p=0.5$) to preserve existing structure.

        \item \textbf{Standard Simulated Annealing:}
        The graph is refined to minimize the unweighted defect sum $E(G) = \text{count}(K_3) + \text{count}(I_7)$. A move consists of flipping a single random edge, accepted via the Metropolis criterion with a geometric cooling schedule ($T_{k+1} = \alpha T_k$).

        \item \textbf{Cost-Driven Termination:}
        The optimization for a specific $n$ runs for a fixed budget (500,000 steps). If the cost reaches 0, the graph is saved as $G_1$ and the search advances to $n+1$. If the cost remains positive, the loop terminates.
    \end{itemize}
\end{summaryboxblue}

\begin{summaryboxblue}[alg:r38]{Search algorithm for {$R(3,8)$}}
    \begin{itemize}
        \item \textbf{Iterative Two-Stage Growth:}
        The algorithm builds the graph vertex-by-vertex ($n \to n+1$). For each new size, it first attempts a lightweight \textbf{Greedy Extension}: the new vertex is connected randomly ($p=0.4$), followed by local edge flips restricted to the new vertex's incident edges to resolve immediate conflicts.

        \item \textbf{Extended Local Search (ELS):}
        If the greedy extension fails to find a valid graph, the algorithm triggers a robust \textbf{ELS Phase}. This involves up to 100 random restarts; in each restart, the algorithm performs \textit{global} edge flips (modifying any edge in the graph, not just the new ones) to escape local minima and reach zero violations.

        \item \textbf{Violation Minimization:}
        The objective function is the unweighted sum of forbidden subgraphs. The search accepts any move (edge flip) that reduces or maintains the violation count, effectively acting as a descent method with plateau traversal:
        $$ E(G) = \text{count}(K_3) + \text{count}(I_8) $$
    \end{itemize}
\end{summaryboxblue}

\begin{summaryboxblue}[alg:r39]{Search algorithm for {$R(3,9)$}}
    \begin{itemize}
        \item \textbf{Algebraic Bootstrap (Circulant Search):}
        The algorithm begins by exhaustively searching for large valid graphs within the space of \textbf{Circulant Graphs}. It generates sum-free sets $S \subset \mathbb{Z}_n$ to define edges $(u, v)$ where $|u-v| \pmod n \in S$. This algebraic structure guarantees $K_3$-freeness by definition and symmetries often yield high Ramsey numbers. 

        \item \textbf{Triangle-Free Growth Strategy:}
        To increase the graph size $n \to n+1$, the algorithm uses a structural heuristic: the new vertex $v_{new}$ is connected exclusively to a \textbf{Maximal Independent Set} $I$ of the current graph. Since no two vertices in $I$ are connected, connecting $v_{new}$ to all $u \in I$ cannot form a triangle.
        $$ N(v_{new}) \leftarrow I \text{ where } I \subseteq V, E(I) = \emptyset $$

        \item \textbf{Constraint Verification:}
        The primary filter is strict $K_3$-freeness. The secondary constraint, $\alpha(G) < 9$, is verified by calculating the independence number (equivalent to the maximum clique in the complement graph $\bar{G}$).

        \item \textbf{Local Refinement:}
        A local search phase attempts to lower the independence number of the graph. It flips edges $(u,v)$ to reduce $\alpha(G)$ but rejects any move that introduces a triangle.
    \end{itemize}
\end{summaryboxblue}

\begin{summaryboxblue}[alg:r44]{Search algorithm for {$R(4,4)$}}
    \begin{itemize}
        \item \textbf{Algebraic Initialization (Paley Graph):}
        The search begins by explicitly constructing the {\bf Paley Graph} of order 17 ($P_{17}$), a known valid graph for $R(4,4)$. Edges are defined based on quadratic residues modulo 17:
        $$ (u, v) \in E \iff (u - v) \pmod{17} \in \{1^2, 2^2, \dots, 8^2\} $$

        \item \textbf{Stochastic Search for Higher Orders:}
        From $n=18$ onwards, the algorithm switches to a stochastic approach. It initializes a random graph $G \sim G(n, 0.5)$ and attempts to refine it into a valid $(4, 4)$-free graph using Simulated Annealing.

        \item \textbf{Symmetric Objective Function:}
        The energy function treats cliques and independent sets symmetrically (since $r=s=4$). The goal is to minimize the total count of 4-cliques and 4-anticliques:
        $$ E(G) = \text{count}(K_4) + \text{count}(I_4) $$

        \item \textbf{Simulated Annealing Loop:}
        The algorithm iteratively flips random edges $(u, v)$. A move is accepted if it lowers the energy $E(G)$ or probabilistically via the Metropolis criterion ($e^{-\Delta E / T}$). If the temperature drops below a threshold $T_{\texttt{min}}$, it resets to $T=1.0$ (reheating) to escape local minima.
    \end{itemize}
\end{summaryboxblue}

\begin{summaryboxblue}[alg:r45]{Search algorithm for {$R(4,5)$}}
    \begin{itemize}
        \item \textbf{Sequential Growth with Multi-Start:}
        The algorithm iterates through target sizes $n$, advancing only after a valid $(4, 5)$-free graph is found. For each size, it launches up to 15 independent restarts to prevent stagnation in local minima.

        \item \textbf{Density-Cycling Initialization:}
        To explore different structural regimes, each restart initializes the graph with a different edge probability, cycling through $p \in \{0.35, 0.40, 0.45\}$. This targets specific density windows where valid $R(4,5)$ graphs are statistically more likely to exist.

        \item \textbf{Standard Simulated Annealing:}
        The refinement process minimizes the unweighted defect sum $E(G) = \text{count}(K_4) + \text{count}(I_5)$. Transitions involve single edge flips, accepted via the Metropolis criterion with a slow geometric cooling schedule ($T_{\texttt{new}} \approx 0.99999 \cdot T_{\texttt{curr}}$).
    \end{itemize}
\end{summaryboxblue}

\begin{summaryboxyellow}[alg:r310]{Search algorithm for {$R(3,10)$}}
    \begin{itemize}
        \item \textbf{Algebraic Bootstrap ($P_{29}$):}
        The search is seeded with a high-quality algebraic construction: the Paley Graph of order 29. This provides a guaranteed symmetric, dense starting point close to the target constraints.
        $$ V = \mathbb{Z}_{29}, \quad (u,v) \in E \iff u-v \in (\mathbb{Z}_{29}^*)^2 $$

        \item \textbf{Iterative Extension Protocol:}
        The algorithm grows the graph incrementally ($n \to n+1$). It only attempts to solve for size $n+1$ after fully resolving the constraints for size $n$ (reaching zero defects).

        \item \textbf{Multi-Start Greedy Initialization:}
        When adding a new vertex $v_{new}$, the algorithm tests three initialization patterns for its incident edges—Empty, Complete, and Random. It applies a greedy local descent to $v_{new}$'s connections to minimize immediate conflicts before starting the full global optimization.

        \item \textbf{Simulated Annealing Refinement:}
        The graph undergoes global Simulated Annealing to minimize the unweighted defect sum $E(G) = \text{count}(K_3) + \text{count}(I_{10})$. Transitions are standard edge flips with a geometric cooling schedule.
    \end{itemize}
\end{summaryboxyellow}

\begin{summaryboxyellow}[alg:r311]{Search algorithm for {$R(3,11)$}}
    \begin{itemize}
        \item \textbf{Hybrid Fractal Initialization:}
        The search is seeded with a diverse pool of starting graphs, including algebraic Paley graphs ($P_{47}$), Cyclic graphs found via sum-free set search, and a unique \textbf{Fractal Construction} that recursively partitions vertices and assigns edge probabilities to create self-similar substructures.

        \item \textbf{Entropic Force Field (Adaptive Scoring):}
        The optimization minimizes a dynamic energy function. The weights $W_{K_3}$ and $W_{I_{11}}$ are not constant; they adjust in real-time based on the ratio of current violations, effectively directing ``force'' against the dominant defect type:
        $$ E(G, t) = W_{K_3}(t) \cdot \text{count}(K_3) + W_{I_{11}}(t) \cdot \text{count}(I_{11}) $$

        \item \textbf{Targeted Edge Flipper (TEF):}
        Instead of proposing purely random edge flips, the algorithm samples a batch of $M$ potential edges. It calculates the weighted $\Delta E$ for all of them and selects the move with the most favorable gradient for the Metropolis acceptance step.

        \item \textbf{Clique Delta Matrix (CDM):}
        To support the high throughput of TEF, the algorithm maintains a dynamic matrix tracking common neighbors, $\mathcal{D}_{ij} = |N(i) \cap N(j)|$. This allows the calculation of $\Delta K_3$ for edge removals in $O(1)$ time.

        \item \textbf{Adaptive Entropic Temperature (AETS):}
        The cooling schedule is driven by a ``frustration index'' (the gap between the current score and the best-known score). This allows the temperature to automatically rise (reheat) when the search stagnates and cool down when rapid progress is being made.
    \end{itemize}
\end{summaryboxyellow}

\begin{summaryboxyellow}[alg:r316]{Search algorithm for {$R(3,16)$}}
    \begin{itemize}
        \item \textbf{Diverse Algebraic Seeding:}
        The search is seeded with candidates from three algebraic constructions:
        \begin{enumerate}
            \item \textit{Standard Cyclic:} Based on sum-free sets in $\mathbb{Z}_n$.
            \item \textit{Product Cyclic:} Based on sum-free sets in $\mathbb{Z}_p \times \mathbb{Z}_q$ (where $n=pq$).
            \item \textit{Multiplicative Coset:} Unions of cosets of multiplicative subgroups in $\mathbb{Z}_n^*$.
        \end{enumerate}

        \item \textbf{Population-Based Extension:}
        A population of valid graphs is maintained, sorted by size. To generate a candidate of size $n+1$, the largest available valid graph is extended by connecting the new vertex to a large independent set of the existing graph (or a random subset thereof).

        \item \textbf{Objective Function:}
        The minimization target penalizes triangles heavily and independent sets linearly:
        $$ E(G) = \text{count}(K_3) + \max(0, \alpha(G) - 15) $$

        \item \textbf{Tabu-Enhanced Local Search:}
        The refinement phase uses Tabu Search to escape local optima. Move selection is strategic:
        \begin{itemize}
            \item If $K_3 > 0$: Prioritize breaking edges involved in triangles.
            \item If $K_3 = 0$ but $\alpha(G) \ge 16$: Prioritize adding edges between vertices in large independent sets.
        \end{itemize}

        \item \textbf{Fast Heuristics:}
        During the intensive local search, $\alpha(G)$ is approximated using a fast randomized bitset heuristic to avoid the expensive exact calculation until the end of the iteration.
    \end{itemize}
\end{summaryboxyellow}

\begin{summaryboxyellow}[alg:r317]{Search algorithm for {$R(3,17)$}}
    \begin{itemize}
        \item \textbf{Heuristic Alpha Filtration:}
        To avoid the computational cost of exact independence number calculations, the algorithm relies on randomized greedy heuristics (using both random and degree-based vertex orderings) to detect independent sets of size $s=17$. Exact verification is reserved for final confirmation.

        \item \textbf{Circulant Bootstrap:}
        The search prioritizes symmetric structures by attempting to construct valid \textbf{Circulant Graphs} generated by sum-free sets $S \subset \mathbb{Z}_n$. This ensures $K_3$-freeness by design, reducing the problem to checking $\alpha(G) < 17$.

        \item \textbf{Smart Extension (Vertex Cloning):}
        When extending a valid graph from $n \to n+1$, the new vertex is not initialized randomly. Instead, it \textbf{clones} the connectivity pattern of a random existing vertex, followed by a slight random perturbation. This preserves the valid structural properties of the parent graph.

        \item \textbf{Targeted Simulated Annealing:}
        The refinement phase minimizes a weighted energy function $E(G) = 100 \cdot \text{count}(K_3) + 10 \cdot \text{count}(I_{17})$. Crucially, the move generator is \textbf{biased}: it selects edges incident to "violating vertices" (those currently part of a clique or large independent set) with 70\% probability, focusing optimization effort where it is most needed.
    \end{itemize}
\end{summaryboxyellow}

\begin{summaryboxyellow}[alg:r319]{Search algorithm for {$R(3,19)$}}
    \begin{itemize}
        \item \textbf{Circulant Bootstrap ($C_{56}$):}
        The search is seeded with a known valid Circulant graph of order 56 ($C_{56}$), defined by a specific set of generators. This provides a massive head-start compared to random initialization.

        \item \textbf{Hybrid Search Strategy:}
        The algorithm alternates between two phases:
        \begin{enumerate}
            \item \textbf{Circulant Search:} Exhaustively searches for larger valid Circulant graphs by testing random generator sets $S \subset \{1, \dots, n/2\}$.
            \item \textbf{Graph Growth:} Greedily extends the best known graph by adding vertices ($n \to n+1$) with connectivity patterns that avoid triangles.
        \end{enumerate}

        \item \textbf{Independent Set Repair:}
        If a candidate graph contains an independent set $I$ of size 19, the algorithm attempts to "repair" it by adding an edge between two vertices $u, v \in I$. This edge destroys the independent set but is only added if it does not create a triangle (i.e., $u, v$ share no common neighbors).

        \item \textbf{Densification:}
        To proactively prevent large independent sets, the algorithm periodically "densifies" the graph by adding random safe edges (edges that do not complete a triangle) until no such edges can be added.
    \end{itemize}
\end{summaryboxyellow}

\begin{summaryboxyellow}[alg:r320]{Search algorithm for {$R(3,20)$}}
    \begin{itemize}
        \item \textbf{High-Quality Initialization ($P_{43}$):}
        The search is jump-started with the Paley graph of order 43 ($P_{43}$), a known $(3, s)$-free graph where $s \ll 20$. This provides a large, dense, valid core structure immediately.

        \item \textbf{Biased Generator Search:}
        The algorithm searches for valid Circulant extensions by generating random generator sets $S$. Crucially, the search is **biased**: new generator sets are often derived from the set $S_{best}$ of the best-known graph, perturbing it slightly rather than starting from scratch.

        \item \textbf{Greedy Alpha Estimation:}
        To filter graphs efficiently, the independence number is estimated using a randomized greedy heuristic. Exact calculation of $\alpha(G)$ is only performed if the greedy estimate is dangerously close to the limit ($18 \le \alpha_{est} < 20$) or if the graph size is small ($n \le 65$).

        \item \textbf{Mutation-Based Growth:}
        Upon finding a new best graph $G_{best}$, the algorithm attempts to grow it to size $n+1$ via mutation. It creates a candidate by cloning $G_{best}$ and adding an isolated vertex, then performs a focused local search (random edge flips and triangle breaking) to resolve new violations.
    \end{itemize}
\end{summaryboxyellow}

\begin{summaryboxyellow}[alg:r321]{Search algorithm for {$R(3,21)$}}
    \begin{itemize}
        \item \textbf{Randomized Local Search (RLS) for $\alpha(G)$:}
        Instead of exact independence number calculation, the algorithm uses a randomized heuristic. It constructs an initial independent set greedily and then performs local vertex swaps (add one vertex, remove its neighbors) to try and increase the set size to 21. If it fails to find an IS of size 21 after many iterations, the graph is potentially valid.

        \item \textbf{Cyclic Graph Construction:}
        The core search strategy focuses on finding large valid \textbf{Cyclic Graphs}. Edges are defined by a generator set $S \subset \{1, \dots, n/2\}$. To ensure $K_3$-freeness, the set $S$ must be sum-free modulo $n$ (i.e., no $x, y, z \in S$ such that $x+y=z$, $x+y+z=0$, etc.).

        \item \textbf{Iterative Generator Augmentation:}
        The sum-free generator set $S$ is built iteratively. The algorithm makes multiple passes over a shuffled pool of candidates $\{1, \dots, n/2\}$. In each pass, it adds any candidate $x$ that does not conflict with the existing members of $S$ (maintaining the sum-free property).

        \item \textbf{Optimized Delta Updates:}
        For local search and graph refinement, the change in the number of triangles ($\Delta K_3$) upon flipping an edge $(u,v)$ is calculated efficiently using vector operations on the adjacency matrix ($\mathbf{A}$) rows:
        $$ \Delta K_3(u,v) = \pm |N(u) \cap N(v)| = \pm \ip{\mathbf{A}_u}{\mathbf{A}_v} $$
    \end{itemize}
\end{summaryboxyellow}

\begin{summaryboxyellow}[alg:r322]{Search algorithm for {$R(3,22)$}}
    \begin{itemize}
        \item \textbf{Cyclic Bootstrap ($C_{101}$):}
        The search bypasses lower-order graphs by initializing with a known valid Circulant graph of order 101 ($C_{101}$) defined by a specific set of generators. This focuses the computational budget on the frontier of the problem space.

        \item \textbf{Sum-Free Generator Construction:}
        The algorithm constructs candidate graphs exclusively within the family of \textbf{Circulant Graphs}. It builds generator sets $S$ by randomly selecting elements $g \in \{1, \dots, n/2\}$ and verifying strict sum-free conditions (e.g., $3g \neq 0$, $2g \notin -S$, $g \notin -(S+S)$). This guarantees $K_3$-freeness algebraically.

        \item \textbf{Randomized Greedy Filtration:}
        To evaluate the independence number $\alpha(G)$, the algorithm employs a fast randomized greedy heuristic. It runs up to 2,000 trials of generating maximal independent sets via random vertex orderings to estimate the lower bound of $\alpha(G)$.

        \item \textbf{Conditional Exact Verification:}
        Exact calculation of the independence number is computationally expensive and is restricted to:
        \begin{enumerate}
            \item Smaller graphs where verification is fast ($n \le 95$).
            \item Promising candidates where the greedy approximation suggests validity ($\alpha_{\texttt{approx}} < 22$).
        \end{enumerate}
    \end{itemize}
\end{summaryboxyellow}

\begin{summaryboxyellow}[alg:r46]{Search algorithm for {$R(4,6)$}}
    \begin{itemize}
        \item \textbf{Cyclic Initialization:}
        The search prioritizes structural regularity by initializing with the best \textbf{Cyclic Graph} found within a short time budget. It generates random sets of "chords" (distances) $D \subset \{1, \dots, n/2\}$ and selects the graph that minimizes the initial defect count.

        \item \textbf{Tabu Search Refinement:}
        The core optimization engine is a \textbf{Tabu Search} that iteratively flips edges to minimize the unweighted sum of violations:
        $$ E(G) = \text{count}(K_4) + \text{count}(I_6) $$
        Recent moves are stored in a Tabu list (tenure $\approx N$) to prevent cycling and encourage exploration.

        \item \textbf{High-Performance Delta Updates:}
        To support high-throughput sampling, the cost change ($\Delta E$) for flipping an edge $(u,v)$ is calculated efficiently using subgraph induction on common neighbors:
        \begin{itemize}
            \item $\Delta K_4$: Depends on edges in the subgraph induced by $N(u) \cap N(v)$.
            \item $\Delta I_6$: Depends on non-edges in the common non-neighborhood of $u, v$.
        \end{itemize}

        \item \textbf{Stagnation Handling (Perturbation):}
        If the search stagnates (no improvement for a dynamic number of iterations based on $N^2$), the algorithm triggers a restart mechanism:
        \begin{itemize}
            \item \textit{High Cost:} Restart with a fresh Cyclic graph.
            \item \textit{Low Cost:} Apply a medium perturbation (flip 25\% of edges).
        \end{itemize}

        \item \textbf{Incremental Growth:}
        The search proceeds sequentially ($n \to n+1$). When extending a valid graph to the next size, the new vertex is connected with a target density of $\approx 0.4$, providing a "warm start" for the Tabu search at the new size.
    \end{itemize}
\end{summaryboxyellow}

\begin{summaryboxyellow}[alg:r47]{Search algorithm for {$R(4,7)$}}
    \begin{itemize}
        \item \textbf{Genetic Circulant Search (Phase 1):}
        The search first explores the structured space of \textbf{Circulant Graphs}. It evolves a population of generator vectors (jump offsets) using a Genetic Algorithm:
        \begin{itemize}
            \item \textit{Seeding:} Initialized with known good generators for $n\in\{35, 36\}$ and random density vectors ($p \approx 0.3-0.5$).
            \item \textit{Operators:} Uses tournament selection, single-point crossover, and bit-flip mutation.
            \item \textit{Elitism:} Preserves the fittest candidate across generations.
        \end{itemize}

        \item \textbf{Targeted Tabu Search (Phase 2):}
        The best candidate from the GA (or an extension of the previous best graph) is refined using Tabu Search. The neighborhood generation is \textbf{focused}:
        $$ \mathcal{N}(G) = \{ (u,v) \mid u,v \in \text{clique}(K_4) \lor u,v \in \text{indep\_set}(I_7) \} \cup \text{Random} $$
        This prioritizes breaking existing violations while maintaining exploration.

        \item \textbf{Perturbation Restart:}
        If the Tabu Search stagnates (no improvement for 100 iterations), the algorithm triggers a "kick" by randomly flipping 10\% of the edges in the current best solution to escape the local minimum.

        \item \textbf{Sequential Extension:}
        If a valid graph for size $n$ is found, the search advances to $n+1$. If scratch construction fails, it attempts to extend $G_{n}$ by adding a new vertex with random connections ($\rho \approx 0.38$) before optimizing.

        \item \textbf{Objective Function:}
        Both phases minimize the unweighted sum of constraints:
        $$ E(G) = \text{count}(K_4) + \text{count}(I_7) $$
    \end{itemize}
\end{summaryboxyellow}

\begin{summaryboxyellow}[alg:r410]{Search algorithm for {$R(4,10)$}}
    \begin{itemize}
        \item \textbf{Circulant Search with Adaptive Probabilities:}
        The primary search mechanism explores \textbf{Circulant Graphs}. The probability of including each generator $g \in \{1, \dots, n/2\}$ is not static; it is learned adaptively. Successful configurations increase the probability of their active generators (reinforcement learning), biasing future searches toward productive subspaces.

        \item \textbf{Near-Miss Repository:}
        The algorithm maintains a repository of "near-miss" graphs (graphs with few violations) for each size $n$. If a perfect graph cannot be found from scratch, the search retrieves a near-miss from the repository (or the best graph from size $n-1$) and attempts to repair it.

        \item \textbf{Hybrid Tabu Repair:}
        The repair process uses a Tabu Search on a hybrid graph representation (Adjacency Matrix + Adjacency List) to enable fast delta updates.
        \begin{itemize}
            \item \textit{Move Selection:} Candidates are chosen by prioritizing edges involved in existing $K_4$ or $I_{10}$ violations, mixed with random edges for exploration.
            \item \textit{Delta Calculation:} $\Delta K_4$ is computed via intersection of adjacency sets; $\Delta I_{10}$ is computed by finding 8-cliques in the common non-neighborhood of the complement graph.
        \end{itemize}

        \item \textbf{Heuristic Filtration:}
        To speed up evaluation, $I_{10}$ checks are performed using a randomized heuristic that estimates the maximum independent set size. Exact verification is only triggered when the heuristic indicates a potential valid solution (cost=0).

        \item \textbf{Sequential Growth:}
        The search advances sequentially. A valid graph $G_n$ is used as a seed for finding $G_{n+1}$, ensuring that the search for larger graphs starts from a high-quality, nearly-valid structure.
    \end{itemize}
\end{summaryboxyellow}

\begin{summaryboxyellow}[alg:r412]{Search algorithm for {$R(4,12)$}}
    \begin{itemize}
        \item \textbf{Algebraic Bootstrap ($P_{127}^*$):}
        The algorithm initializes with a high-quality algebraic graph of order 127 based on cubic residues in $\mathbb{F}_{127}$. Specifically, $V = \mathbb{F}_{127}$ and $(u,v) \in E \iff (v-u) \in \{x^3 \mid x \in \mathbb{F}_{127}^*\}$. This graph is known to be $(4, 12)$-free.

        \item \textbf{Extension Search (Phase 1):}
        To reach $n=128$, the algorithm searches for a valid "extension mask" (neighborhood vector for the new vertex $v_{128}$). This is formulated as minimizing a cost function:
        $$ E(m) = \text{count}(K_3 \in N(v_{128})) + \text{count}(I_{11} \in V \setminus N(v_{128})) $$
        A short Simulated Annealing phase optimizes this mask to find a low-defect starting configuration.

        \item \textbf{Targeted Defect Repair (Phase 2):}
        The 128-vertex graph is refined using a specialized Simulated Annealing process.
        \begin{itemize}
            \item \textbf{Defect Tracking:} The algorithm periodically enumerates all $K_4$ and $I_{12}$ subgraphs and maps them to their constituent edges.
            \item \textbf{Biased Sampling:} Edge flip candidates are not chosen uniformly; edges participating in more defects are weighted higher for selection.
            \item \textbf{Fast Delta Update:} The cost change $\Delta E$ is computed efficiently by inspecting only the local neighborhood of the flipped edge (checking for $K_2$ in common neighbors or $I_{10}$ in common non-neighbors).
        \end{itemize}
    \end{itemize}
\end{summaryboxyellow}

\begin{summaryboxyellow}[alg:r417]{Search algorithm for {$\ramsey(4,17)$} lower bound}
\begin{itemize}
\item \textbf{Multiplier-Based Orbit Construction:} The search is restricted to highly symmetric, vertex-transitive circulant graphs of size $n$, fully defined by a difference set $S \subset \{1, \dots, \lfloor n/2 \rfloor\}$. To drastically compress the search space, the difference set is constructed as a union of symmetric orbits generated by modular multiplication with a coprime multiplier $a$, allowing the algorithm to toggle entire symmetry groups simultaneously.
\item \textbf{Vertex-Transitive Reduction:} Because the graph is vertex-transitive, global structural constraints are reduced to localized checks strictly around vertex 0. The presence of forbidden subgraphs is mapped directly to the neighborhood and non-neighborhood of this single vertex:
\begin{align*}
K_4 \text{ freeness} &\iff \text{No } K_3 \text{ in the subgraph induced by the neighborhood } S \\
I_{17} \text{ freeness} &\iff \text{No } I_{16} \text{ in the subgraph induced by the non-neighborhood } \mathbb{Z}_n \setminus (S \cup \{0\})
\end{align*}
\item \textbf{Scaled Seeding and Multiplier Tournament:} As the algorithm sequentially tests graph sizes $n \in [199, 215]$, it runs a "tournament" evaluating various coprime multipliers $a$ to find the most favorable orbit structures. To accelerate convergence for each new $n$, the initial search state is intelligently seeded by projecting and scaling the successful difference set from the previously solved graph size onto the new orbit layout.
\item \textbf{Simulated Annealing with External Acceleration:} The core search engine uses simulated annealing to toggle orbits in or out of $S$, minimizing an energy function based on $K_4$ and $I_{17}$ violations. Because finding large independent sets is computationally expensive, the localized subgraph checks are dynamically offloaded to a highly optimized external C++ library (`cliques\_wrapper`), enabling the rapid evaluation necessary for the annealing loops.
\end{itemize}
\end{summaryboxyellow}

\begin{summaryboxyellow}[alg:r57]{Search algorithm for {$R(5,7)$}}
    \begin{itemize}
        \item \textbf{Vertex-Transitive Search Space:}
        The algorithm restricts the search to \textbf{Cayley Graphs} based on cyclic groups. Because these graphs are vertex-transitive, the structural properties of every vertex are identical. This allows the cost function to be evaluated by inspecting only the local neighborhood of a single vertex (vertex 0), massively reducing computational overhead.

        \item \textbf{Localized Cost Function:}
        Instead of counting all cliques in the graph, the algorithm minimizes a localized proxy score:
        \begin{itemize}
            \item $K_5$: Approximated by finding $K_4$ in the neighborhood of vertex 0.
            \item $I_7$: Approximated by finding $I_6$ in the non-neighborhood of vertex 0.
        \end{itemize}
        $$ E(S) \approx \frac{n}{5} \cdot \text{count}(K_4 \in N(0)) + \frac{n}{7} \cdot \text{count}(I_6 \in V \setminus N(0)) $$

        \item \textbf{Stochastic Generator Optimization:}
        The generator set $S \subset \{1, \dots, n/2\}$ is optimized using Simulated Annealing. Moves consist of adding or removing a generator from $S$.

        \item \textbf{Adaptive Local Refinement:}
        If a graph with low (but non-zero) defects is found, it undergoes a refinement phase using \textbf{adaptive edge weighting}. Edges that frequently participate in beneficial flips are assigned higher weights (interest), guiding the local search toward critical structures.
    \end{itemize}
\end{summaryboxyellow}

\begin{summaryboxyellow}[alg:r58]{Search algorithm for {$R(5,8)$}}
    \begin{itemize}
        \item \textbf{Vertex-Transitive Reduction:}
        The search relies heavily on \textbf{Circulant Graphs}. Exploiting vertex transitivity, the global clique count is reduced to local subgraph counting around vertex 0:
        \begin{itemize}
            \item $K_5$ count $\rightarrow \frac{n}{5} \times \text{count}(K_4 \text{ in } N(0))$
            \item $I_8$ count $\rightarrow \frac{n}{8} \times \text{count}(K_7 \text{ in } \overline{V \setminus N(0)})$
        \end{itemize}

        \item \textbf{Fractal Mirror Initialization:}
        When extending a valid graph from $n \to n+1$, the algorithm uses a "Fractal Self-Similarity" heuristic. Instead of random initialization, the new vertex $v_{new}$ mirrors the connectivity of a randomly selected existing vertex $v_{mirror}$, with a 15\% noise factor. This preserves the local Ramsey-safe substructures of the parent graph.

        \item \textbf{Parallelized Simulated Annealing:}
        The solution space is explored concurrently using a multiprocessing pool. Each worker runs an independent Simulated Annealing chain, swapping generator offsets to minimize defects.

        \item \textbf{Delta-Update Local Search:}
        For general (non-circulant) refinement, edge flips are evaluated using efficient local counting:
        \begin{itemize}
            \item $\Delta K_5$: Calculated by finding $K_3$ in the common neighborhood.
            \item $\Delta I_8$: Calculated by finding $I_6$ in the common non-neighborhood.
        \end{itemize}

        \item \textbf{Sequential Hybrid Strategy:}
        The workflow alternates: find a valid base using the fast Circulant search, then push the size limit using the general Local Search with the fractal extension strategy.
    \end{itemize}
\end{summaryboxyellow}

\begin{summaryboxyellow}[alg:r59]{Search algorithm for {$R(5,9)$}}
    \begin{itemize}
        \item \textbf{Circulant Subspace Search:}
        The algorithm exclusively searches the space of \textbf{Circulant Graphs}. The search space is defined by the binary vector of generators $\{g_1, \dots, g_{\lfloor n/2 \rfloor}\}$. This algebraic structure is chosen because known Ramsey graphs frequently exhibit high symmetry.

        \item \textbf{Localized Scoring:}
        The objective function exploits vertex transitivity to reduce the cost of calculating global defects. Instead of full enumeration, it checks only the neighborhood of vertex 0:
        \begin{itemize}
            \item $K_5$: Approximated by finding $K_4$ in the neighborhood $N(0)$.
            \item $I_9$: Approximated by finding $I_8$ (i.e., $K_8$ in complement) in the non-neighborhood $\overline{N(0)}$.
        \end{itemize}

        \item \textbf{Iterated Local Search (ILS):}
        The generator vector is optimized using a two-tier strategy:
        \begin{enumerate}
            \item \textbf{Tabu Search (Inner Loop):} Greedily flips generators to minimize the local score, using a Tabu list to prevent cycling. Aspiration criteria allow tabu moves if they yield a new global best.
            \item \textbf{Perturbation (Outer Loop):} When Tabu Search settles, the solution is "kicked" by randomly flipping a small fraction (2-10
        \end{enumerate}

        \item \textbf{Aggressive Growth Strategy:}
        The search iterates graph sizes $n$. Upon finding a valid graph at size $n$, it attempts an aggressive "jump" (e.g., $n \to n+8$) to rapidly probe upper bounds. If the jump fails, it falls back to incremental growth ($n \to n+1$) to fill the gap.
    \end{itemize}
\end{summaryboxyellow}

\begin{summaryboxyellow}[alg:r67]{Search algorithm for {$R(6,7)$}}
    \begin{itemize}
        \item \textbf{Quadratic Residue Initialization ($P_p$):}
        The search is seeded with a \textbf{Paley Graph} constructed from quadratic residues modulo a prime $p$. This algebraic structure ensures high symmetry and avoids small cliques by definition ($K_k$-free for $k \approx \ln p$).
        $$ (u, v) \in E \iff (v - u) \pmod p \in \{x^2 \mid x \in \mathbb{F}_p^*\} $$

        \item \textbf{Asynchronous Tabu Search:}
        The graph is refined using Tabu Search. To handle the high computational cost of counting $K_6$ and $I_7$, the algorithm employs \textbf{Asynchronous Re-evaluation}: full violation lists are recomputed only periodically ($T \approx n/3$), while inner loops rely on cached (potentially stale) data to maximize iteration speed.

        \item \textbf{Genetic Column Generation:}
        
        When extending the graph ($n \to n+1$), the connectivity of the new vertex is optimized using a \textbf{Genetic Algorithm}. The population consists of binary "connection vectors" evolved via crossover and mutation, initialized based on the connectivity patterns of existing vertices (mimicry) rather than pure randomness.

        \item \textbf{Chaotic Descent Restart:}
        If the local search stagnates, a "Chaotic Descent" mechanism is triggered. It perturbs the graph using a deterministic but distinct quadratic residue mapping based on a different prime $q$, effectively scrambling the adjacency matrix while retaining some quasi-random algebraic properties.

        \item \textbf{Violation-Frequency Sampling:}
        Edge flip candidates are selected based on a frequency map: edges that appear most frequently in the set of all current $K_6$ and $I_7$ subgraphs are prioritized for flipping.
    \end{itemize}
\end{summaryboxyellow}

\end{document}